\PassOptionsToPackage{rgb}{xcolor}
\documentclass{article}
\usepackage[
journal=JST,
lang=american
]{ems-journal}

\usepackage{mathbbol}
\usepackage{verbatim}

\newtheorem{defi}{Definition}[section]
\newtheorem{thm}{Theorem}[section]

\newtheorem{rem}{Remark}[section]
\newtheorem{rems}{Remarks}[section]
 
\newtheorem{lemma}{Lemma}[section]
\newtheorem{cor}{Corollary}[section]

\newcommand{\extd}{\delta}
\newcommand{\secmat}{\mathfrak{A}}

\numberwithin{equation}{section}

\title
{On topological
bound
states and secular equations for quantum-graph eigenvalues}
\titlemark{Topological quantum-graph states}

\emsauthor{1}{Evans M. Harrell II}{E.~Harrell}
\emsaffil{1}{School of Mathematics,
Georgia Institute of Technology,\
Atlanta GA 30332-0160, USA. \email{eharrell@duck.com}}

\emsauthor{2}{Anna V. Maltsev}{A.~Maltsev}
\emsaffil{2}{School of Mathematical Sciences,
Queen Mary University of London,\
London E1 4NS, UK. \email{a.maltsev@qmul.ac.uk}}

\begin{document}

\maketitle



\section*{Abstract}
Quantum graphs without interaction which contain equilateral cycles possess
``topological'' bound states which
do not correspond to zeroes of one of the two variants of the secular equation for quantum graphs.  Instead, their eigenvalues lie in the set of singularities of the vertex-scattering secular matrix.
This observation turns out to be representative of a wider phenomenon.
We introduce a notion of topological bound states and show that they are linear combinations of functions supported on
generators of the fundamental group of the graph (hence the ``topological" in the name), including for graphs that have interactions on the edges.
Using an Ihara-style theorem, we elucidate
the role of such topological bound states in the spectral analysis of quantum graph Hamiltonians using secular matrices.
En route we determine the set of
the fixed vectors of the bond-scattering matrix.

This work is dedicated to E.B. Davies on the occasion of his 80th birthday and in honor of his
important contributions to the theory of quantum graphs, e.g., \cite{DaExLi,Da13} and of his
broad and influential
work on spectral theory, e.g., \cite{Da89,Da95}.

\section{Introduction}

With the rise of big data, graph theory in the guise of network science has found interesting
new applications in fields of study as diverse as biology, physics, and sociology.
Some applications in microelectronics and neural networks, known as quantum graphs, require edges
not only to connect vertices but also to carry solutions of differential equations.
Reference works about the theory of quantum graphs include \cite{BeKu,Bor,Kur}.

Many fundamental questions about networks relate to the spectral properties of operators acting on them, such as the
discrete Laplacian.
For example,
It is known as part of Hodge theory for combinatorial
graphs that eigenfunctions of discrete Laplacians often localize on the
topological ``holes''; see e.g. \cite{BaSa,BeLo,bianconi2021topological, Lim20}
and references therein for this phenomenon and some consequences.
The possibility of eigenfunctions that localize on compact subsets such as cycles is likewise well-known in the theory of
quantum graphs.
In a
sense that has been made precise in
\cite{BeLi,Alon}, quantum-graph eigenfunctions ``generically'' do not vanish on vertices and consequently are not compactly supported,
yet compactly
supported eigenfunctions
occur rather often in commonly considered models.
Much of the literature referencing such states of quantum graphs relates them to
scattering resonances, following the work of Exner and Lipovsk\'y \cite{ExLi}, who suggested that
near-real ``topological' could grow out of compactly supported bound states; this idea has  been further explored in
\cite{CdV,GSS,CdVTru}.  Compactly supported states have also been related to the issue of ``scarring'' in \cite{KotSch}.

The main result of this article is to identify a topological mechanism that gives rise to compactly supported eigenfunctions on quantum graphs.  After some background in this Introduction we
review the use of secular equations involving the ``bond scattering matrix'' $\mathcal{S}$
as a tool for spectral analysis on quantum graphs, innovating in \eqref{SecEq} and the following discussion in
order to incorporate general interactions on the edges.  The operator $\mathcal{S}$ is a purely topological quantity unrelated to the eigenparameter.   In \S
\ref{s: BSM} we present a useful factorization of $\mathcal{S}$ and several other quantities of interest in Theorem \ref{Facto},
and we then characterize the fixed points of $\mathcal{S}$ in Theorem \ref{FixedSvecs}.
In particular, we prove that, other than the constant vector, the
fixed points correspond to discrete harmonic vector fields and
have bases supported on ``holes'' that generate the fundamental group of the combinatorial aspect of the metric graph.
Section \ref{s: TQGS} turns to quantum graphs, where something happens on edges in addition to purely topological relations.
We identify a topological condition related to the eigenfunctions that fail to be identified in one of the versions of the
secular equation, and use the condition
as the basis of Definition \ref{TQGSDeff}.  Theorem \ref{t:topoCharacterization} precisely characterizes topological quantum-graph states and shows that they are
supported on holes for reasons similar to the case of the fixed vectors of $\mathcal{S}$.
Theorem \ref{t:classification} accounts for the different types of quantum-graph eigenvalues and sets up a vertex secular equation that includes
general edge interactions; an expicit form is worked out in the following corollary.  We close with a section consisting of several instructive case studies.

Heretofore topological bound states appear to be
systematically discussed only in symmetric cases with no interaction on the edges.
Analysis of compactly supported states on metric graphs began, at least implicitly, with the 1985 work of von Below \cite{Below},
in which
 it was shown that the spectral analysis of equilateral metric graphs, with no quantum interaction and standard Kirchhoff vertex conditions, can be reduced by factorization to the analysis of the underlying combinatorial graph.  The equilateral situation has been further analyzed in greater generality
 in \cite{BeMu,Pank,BBJL}.
The equilateral case is accessible to analysis because symmetry can be used to reduce the spectral analysis to the study of Sturm-Liouville
problems on individual edges with Dirichlet or respectively Neumann conditions at the vertices.  The noninteracting case factorizes so that the spectrum of the Laplacian on the metric graph is completely determined by that of the combinatorial Laplacian, for which compactly supported eigenvectors are typical.

As an instructive case consider
the equilateral tetrahedral metric graph with edge length $\pi$, no quantum interaction, and standard Kirchhoff vertex conditions,
i.e., eigenfunctions
$\psi_\ell$ are continuous at the vertices, and the sum of their outgoing derivatives at each vertex is $0$.  On the
edges, $- \psi_\ell^{\prime\prime} = k^2 \psi_\ell$.   An exercise exploiting symmetry
reveals that the eigenfunctions may be divided into those that are symmetric with respect to permutation of the edges incident to each vertex, and those that vanish at each vertex, and thus the eigenvalues equal those of the interval with respectively Neumann and Dirichlet boundary conditions, i.e., $k$ takes on integer values.
Furthermore, for the category of eigenfunctions that vanish at vertices it turns out to be possible to find eigenbases that consist of functions supported on quadrilateral cycles of the tetrahedron.  (When $k$ is an even positive integer
there are eigenbases consisting of functions supported on triangular faces; see Case Study \ref{Tet}, below.)
As to the multiplicity, the dimension of these eigenspaces is 3, which equals the
first Betti number $\beta_1 : = |E| - |V| + 1 = 6 - 4 + 1$ of the discrete tetrahedral graph, i.e., the
 number of generators of its fundamental group.  As we show below in Theorem \ref{FixedSvecs}, this association with the fundamental group is not a coincidence.

 Are compactly supported states like those of the tetrahedron merely exceptional
 cases due to high symmetry, or are they instances of a general phenomenon?
By embedding such highly symmetric examples in larger, quite arbitrary graphs attached at vertices where eigenfunctions vanish, it is easy to see that the possibility of eigenfunctions
supported on the graph's ``holes'' is by no means restricted to equilateral graphs.
One of our purposes in this article is to explore conditions for the existence of compactly supported eigenstates that are associated with fundamental group of the combinatorial graph on which the quantum graph is built, and to characterize them.

A second purpose of this article is to clarify some aspects of the usual method of determining quantum-graph eigenvalues via a secular determinant \cite{BeKu}, which dates from the pioneering work of Kottos and Smilansky \cite{KS99}.  As already realized in \cite{KS99}, there are (at least) two distinct ways to set up a secular matrix, one focusing on ``vertex scattering'' and one focusing on ``bond scattering''  (The latter has become the more widespread approach, e.g., see \cite{BeKu}; The reader is referred to \cite{KS99,BeKu} for the motivation and
usual derivation of the secular matrices.)  As part of our analysis we build an alternative
algebraic bridge between these two secular matrices following a technique of Ihara \cite{Iha, Sun08},
in the course of which we show that the two secular equations are not completely equivalent.  Indeed, although the vertex-scattering version of the secular matrix is smaller and therefore in some ways more efficient than the bond-scattering version, it does not account for certain topological quantum-graph eigenstates analogous to the tetrahedral states supported on cycles.

The great majority of articles on quantum-graph spectral problems using a bond-scattering analysis look only at the free Schr\"odinger  equation
$-u'' = k^2 u$ on the edges, or include a purely magnetic interaction, which in
one-dimension can be locally gauged into a free Schr\"odinger form with a change of variable.
Here
we take care to spell out how to incorporate interactions
on the edges by allowing
a general transfer matrix $\Phi_{\bf e}$.
We imagine $\Phi_{\bf e}$
as the solution operator for a second-order ordinary differential equation
$-u'' + V(x) u = k^2 u$, but it is far from restricted to this particular form.

We recall that a quantum-graph eigenfunction $\psi(x; k)$ is an $L^2$-normalized continuous
function that satisfies
\begin{equation}\label{EVE}
\left(- \frac{d^2}{dx^2} + V(x)\right) \psi(x; k) = k^2 \psi(x; k)
\end{equation}
on the edges of a metric graph $\Gamma$,
and certain conditions at the vertices.
In this article we confine ourselves to
Kirchhoff (a.k.a.\ Neumann-Kirchhoff \cite{Ber} or ``standard'') vertex conditions,
according to which
functions are continuous at vertices and
the sum of their outgoing derivatives at each vertex is $0$. We refer to \cite{BeKu,Ber,Post} for background and precise definitions of these operators.

In \cite{KS99} the secular determinant was derived by regarding the vertex conditions and transfer operators on bonds as a linear system and using the methods of linear algebra to reduce them to a convenient form in which the combinatorial part is captured in a $2m \times 2m$ matrix acting on the set of oriented edges, known as the {\em bond-scattering matrix} $\mathcal{S}$.
(Here and henceforth we set $m := |E|$, and
$n:=|V|$.)
\begin{defi}\label{d:bond-scattering}
The bond-scattering matrix $\mathcal{S}$ is given as follows:
\begin{itemize}
\item $\mathcal{S}_{{\bf e^\prime}{\bf e}} = 0$ unless the terminal vertex  $t({\bf e})$ equals the source vertex $s({\bf e^\prime})$.
\item Otherwise, if ${\bf e^\prime} \ne - {\bf e}$, where $- {\bf e}$ designates the reversal of $- {\bf e}$, then $\mathcal{S}_{{\bf e^\prime}{\bf e}} = \frac{2}{d_{t({\bf e})}}$.
\item
If ${\bf e^\prime} = - {\bf e}$, then
 $\mathcal{S}_{-{\bf e}{\bf e}} = \frac{2}{d_{t({\bf e})}} - 1$.
 \end{itemize}
\end{defi}
We observe that the matrix $\mathcal{S}$ is unitary and doubly stochastic.
(See the discussion of $\mathcal{S}$ in \cite{BeKu} for further details and
\cite{BeWi} for some additional spectral properties of $\mathcal{S}$.)
The bond-scattering secular determinant for the eigenvalues of a quantum graph joins the topological aspects of the graph with the analysis on the edges through a condition of the form
\begin{equation}\label{SecEq1}
\zeta_S(k) := \det (\mathbb{1} - \Phi(k) \mathcal{S}) = 0,
\end{equation}
in which $\Phi$ depends analytically on the effective eigenparameter $k$ and makes no reference to the topology of the graph, whereas the
$k$-independent {\em bond-scattering matrix} $\mathcal{S}$ is purely combinatorial. Solutions of equation \eqref{SecEq1} are the eigenvalues of the quantum graph \cite{KS99, Ber}.

Allowing a nontrivial
scalar interaction on the edges brings a small additional complication.  Without scalar interactions Kottos and Smilansky in \cite{KS99} were able
to choose the matrix $\Phi$ in \eqref{SecEq1} in the form $\textrm{diag}(\exp(i k \ell_{\bf e}))$, where $\ell_{\bf e} = \ell_{- \bf e}$
 is the length of the oriented edge ${\bf e}$ and  $k^2$ is an eigenvalue
of the quantum graph when the secular determinant vanishes.  This, however, relies on a simplification that is special to $-\psi'' = k^2 \psi$,
with which the $2 \times 2$ fundamental solution matrix can be encoded as a single complex-valued function
$\exp(i k \ell_{\bf e})$.  (See Remarks \ref{rems:Phi} (\ref{rem:identifyi}) below.)

In order to handle  ODEs with interactions such as  $-\psi'' + V(x) \psi = k^2 \psi$,
we replace
Eq. \eqref{SecEq1} by
\begin{equation}\label{SecEq}
\det (\mathbb{1}_{4m} - \Phi(k)\hat{\mathcal{S}}) = 0,
\end{equation}
where $\hat{\mathcal{S}} := \mathcal{S} \otimes \mathbb{1}_{2}$
and $\Phi(k)$ is a block-diagonal matrix composed of $2 m$  $2 \times 2$ matrices $\Phi_{\bf e}(k)$.
(When the details about the dimensions of the matrices are unimportant
we may abbreviate \eqref{SecEq} in the form \eqref{SecEq1} for simplicity.)
Since block-diagonal matrices will occur frequently, we choose a labeling convention and introduce the notation

 \[
 \textrm{BMat} [M_{\bf e}] :=
 \begin{bmatrix}
M_{{\bf e}_1} &  &  & &&&&\\
   &M_{{\bf e}_2} &  & &&&&\\
   &  &  \ddots & &&&&\\
   &  &   & M_{{\bf e}_m}   &&&&\\
   &&&& M_{-{\bf e}_1}  &&&\\
     &&&& &M_{-{\bf e}_2} &&\\
        & & &&& &  \ddots &\\
        & &&&& &&M_{-{\bf e}_m}
 \end{bmatrix}.
 \]

Recall that a fundamental solution set for $-\psi'' + V(x) \psi = k^2 \psi$ consists of solutions $\psi_{1}, \psi_2$ such that at the
beginning end of the edge ${\bf e}$ parametrized as $x=0$, $\psi_{1}(0) = 1$, $\psi_{1}'(0) = 0$, while $\psi_{1}(0) = 0$, $\psi_{1}'(0) = 1$.
The general solution matrix, or {\em transfer matrix}  $\Phi_{\bf e}(k)$ on the edge ${\bf e}$ takes on the form

\begin{equation}\label{e:genSol}
\Phi_{\bf e}(k) = \begin{pmatrix}
\psi_1(k,\ell_{\bf e}) & \psi_2(k,\ell_{\bf e}) \\
\psi_1^\prime(k,\ell_{\bf e}) & \psi_2^\prime(k,\ell_{\bf e})
\end{pmatrix}
\end{equation}
at the terminal vertex of ${\bf e}$. Canonically, with $V=0$, this convention gives $U$ the form
\begin{equation}\label{e:V0Sol}
\Phi_{\bf e}(k) = \begin{pmatrix}
\cos(k \ell_{\bf e}) & \frac{1}{k} \sin(k \ell_{\bf e})  \\
-k \sin(k \ell_{\bf e})  & \cos(k \ell_{\bf e})
\end{pmatrix}.
\end{equation}

In what follows, we assume only
the following conditions on the transfer matrices  $\Phi(k)$:

\noindent
{\bf Assumptions}

\begin{enumerate}[A.]\label{PhiProps}
\item
In a basis of oriented edges ${\bf e}$, $\Phi(k)$ is a block-diagonal matrix consisting of $2 \times 2$
invertible matrices
$\Phi_{\bf e}(k)$.
\item
For standardization, the first entry of the vector space ${\mathbb C}^2$ attached to the directed edge ${\bf e}$ corresponds to the value of a quantum-graph eigenfunction at the source vertex $s({\bf e})$ and the second entry to its outgoing derivative at $s({\bf e})$.
\item
Each block $\Phi_{\bf e}(k)$ is an entire matrix-valued function of $k$.
\end{enumerate}

Quantum states $\psi$ are functions in $L^2$ of the metric graph,
but they are determined by
vectors $x$ in the null space of the matrix
appearing in the determinant in \eqref{SecEq}.  We therefore make the dollowing definition.

 \begin{defi}   If $\psi$ is a quantum-graph eigenfuction, then any non-zero vector $x$ in the nullspace of the matrix
 $4m \times 4m$ matrix $\mathbb{1}_{4m} - \Phi(k)\hat{\mathcal{S}}$ will be referred to as a
 {\em counterpart vector} to $\psi$, cf.,  \eqref{SecEq}.
Thus the counterpart vector is a family of
2-component vectors
attached to each directed edge ${\bf e}$, for which $x_1(s({\bf e}) = \psi$ and $x_2(s({\bf e}) = \psi^\prime(s({\bf e})$.  This ensures that
such an $x$ determines the corresponding $\psi$ up to normalization.
 \end{defi}

\begin{rems}\label{rems:Phi} \,

\begin{enumerate}

\item
When the first entry of a $2$-vector attached to ${\bf e}$ vanishes, we say
that the corresponding quantum state $\psi$
satisfies a {\em Dirichlet} condition at $s({\bf e})$.
Analogously, if the second entry of a counterpart $x$ vanishes it is equivalent to a
{\em Neumann} boundary condition at $s({\bf e})$.  Periodic boundary conditions will also arise naturally; see the remark
following \eqref{InverseCalc}.
\item
When setting up the solution matrices on the edges $\Phi_{\bf e}(k)$,
our convention is to preserve the orientation on ${\bf e}$, i.e., the entries of $\Phi_{\bf e}(k) x$ correspond
respectively to the value of the eigenfunction and its {\em incoming}
rather than outgoing derivative at $t({\bf e})$.
Because the orientation of $-{\bf e}$ is the reverse of that of ${\bf e}$, the $x$-derivative on $-{\bf e}$ is the negative of
the $x$-derivative on ${\bf e}$.
Consequently the $2 \times 2$ blocks $\Phi_{\bf e}(k)$ satisfy:
\begin{equation}\label{OrRev}
\Phi_{- \bf e}(k) =
\begin{pmatrix}
1 & 0\\
0 & -1
\end{pmatrix}
\Phi_{\bf e}^{-1}(k)
\begin{pmatrix}
1 & 0\\
0 & -1
\end{pmatrix}.
\end{equation}

\item\label{rem:identifyi}
If
$
\Phi_{\bf e}(k) = \small\small
\begin{pmatrix}
\cos k x & \frac{\sin k x}{k}\\
-k \sin k x & \cos k x
\end{pmatrix}
$
then
\begin{equation}\label{e:reversal}
\Phi_{- \bf e}(k) = \Phi_{\bf e}(k).
\end{equation} Furthermore, when $V=0$ we have
\begin{equation}
\Phi_{\bf e}(k) = \mathbb{1}_2\cos k + \begin{pmatrix}
0 & \frac{1}{k}\\
-k  & 0
\end{pmatrix}
\sin k.
\end{equation}
As $\det \begin{pmatrix}
0 & \frac{1}{k}\\
-k  & 0
\end{pmatrix} = 1$ and $\begin{pmatrix}
0 & \frac{1}{k}\\
-k  & 0
\end{pmatrix}^2 = -\mathbb{1}_2$ we can identify
\begin{equation}\label{e:i}
i = \begin{pmatrix}
0 & \frac{1}{k}\\
-k  & 0
\end{pmatrix}
\end{equation}
where $i$ is the imaginary unit. This allows us to reduce the secular matrix to the $2m \times 2m$ form in  \cite{Ber} rather than the $4m \times 4m$ matrix that we introduced here to account for a possibility of a non-zero $V$.  In fact \eqref{e:reversal} remains true with a scalar potential energy that is symmetric on ${\bf e}$, but
it is not true in general.
\end{enumerate}

\end{rems}

\section{Structure and fixed vectors of the bond-scattering matrix}\label{s: BSM}

The bond-scattering matrix given in Definition \ref{d:bond-scattering} and used to find quantum graph eigenvalues via \eqref{SecEq1} is related to the other matrices that arise in algebraic graph theory including the adjacency matrix, the graph Laplacian, the non-backtracking matrix, and the discrete heat kernel; we detail some of these connections below in Theorem \ref{Facto}.
Because of complications arising from a reliance on oriented edges, we shall write such matrices in terms of the
{\em oriented} incidence matrix of the underlying combinatorial graph.

In discrete differential geometry, the incidence matrix is interpreted as the exterior derivative and in cohomology theory as the coboundary operator;
for discussions of these matters from different points of view we refer to \cite{BeLo,Sun08,Lim20}.
A version of the incidence matrix adapted to the space of oriented edges can be defined as follows:
Let
 $\extd^*_t$ be the $n \times 2m$ matrix associating a given edge with its terminal vertex $t({\bf e})$,
 in that  $(\extd^*_t)_{t({\bf e}),{\bf e}} = 1$ with all other entries $0$.
Correspondingly let  $\extd^*_s$ associate ${\bf e}$ with its source vertex $s({\bf e})$.  (The difference $\extd^*_t - \extd^*_s$ then gives the full oriented incidence matrix,
 and the usual discrete derivative is its
adjoint, $\extd_t - \extd_s$.  Terminology in the literature is not universal, conventions varying as to sign and whether the incidence matrix is represented as an $n \times 2m$ or an $n \times m$ matrix.)   Denote the edge-reversal operator by $\tau: {\bf e}  \to -{\bf e}$.  Note that if the oriented edges are ordered so that the second $m$ of them are the reversals of the first $m$ in the same order, then $\tau$ has the block form
 \[
  \begin{pmatrix}
   0 & \mathbb{1}_m\\
  \mathbb{1}_m & 0
  \end{pmatrix}.
\]
From the definitions, $\extd^*_t \tau = \extd^*_s$ and $\extd^*_s \tau = \extd^*_t$.


\begin{thm}\label{Facto}
The bond-scattering matrix can be written as:
\begin{equation}\label{Sdef}
\mathcal{S} = 2\extd_s \mathcal{D}^{-1}\extd^*_t - \tau,
\end{equation}
where $\mathcal{D} = \extd^*_s \extd_s = \extd^*_t \extd_t$ is the diagonal matrix of the degrees $d_v$ of the vertices.
\begin{table}[H]
\caption{Some additional matrix factorizations}\label{FactTable}
\begin{center}
\begin{small}
\begin{tabular}{| l | l | l | l | }
\hline
{\bf matrix} &  {\bf formula} & {\bf dimensions} & {\bf notes}\\
\hline
Incidence  &  $\extd_t^* - \extd_s^*$ & $n \times 2m$ & oriented\\
Adjacency  &  $\mathcal{A} = \extd^*_t \extd_s$ & $n \times n$ &\\
Degrees   &   $\mathcal{D} = \extd^*_t \extd_t = \extd^*_s \extd_s$ & $n \times n$ & diag$(d_v)$\\
Discrete Laplacian   & $\mathcal{L} = \mathcal{D} - A = \extd^*_t (\extd_t - \extd_s) = (\extd_s^* - \extd_t^*) \extd_s$ &$n \times n$ & $\mathcal{L} \ge 0$\\
Normalized Laplacian & $\mathcal{D}^{-1/2}   (\extd_t^* - \extd_s^*) (\extd_t - \extd_s)\mathcal{D}^{-1/2}$ & $n \times n$ & \\
Laplacian on forms  &$(\extd_t -\extd_s) \mathcal{D}^{-1}(\extd^*_t - \extd^*_s )$& $2m \times 2m$  & reduces to $m \times m$\\
=1-down Laplacian &&&\\
Edge reversal & $\tau =
  \begin{pmatrix}   0 & \mathbb{1}\\  \mathbb{1} & 0  \end{pmatrix}$ & $2m \times 2m$

  & (or $4m \times 4m$) \\
Bond-scattering & $\mathcal{S} = 2\extd_s \mathcal{D}^{-1}\extd^*_t - \tau$ & $2m \times 2m$ & unitary\\
Nonbacktracking & $\mathcal{H} =  \extd_t  \extd^*_s$ & $2m \times 2m$ & self-loops are excluded\\
\hline
Bond transfer & $\Phi(k) = \textrm{BMat}(\Phi_{\bf e}(k))$ & $4m \times 4m$ & edge solution matrices\\
Secular & $ \mathbb{1}_{4m} - ( \mathcal{S} \otimes \mathbb{1}_{2}) U(k) $ & $4m \times 4m$ & edge version\\
\hline
\end{tabular}
\end{small}
\end{center}
\end{table}
\noindent

Basic facts about the edge-reversal operator $\tau$ include:

$\tau^2 = \mathbb{1}$; $\extd_t^* \tau = \extd_s^* $  and $\extd_s^* \tau = \extd_t^* $;
$\tau \mathcal{S} \tau = \mathcal{S}^{-1} = \mathcal{S}^*$

\end{thm}

\begin{proof}
Each of these identities is a straightforward exercise.
\end{proof}

The next theorem characterizes the fixed points of the operator $\mathcal{S}$ as harmonic $1$-forms on the combinatorial graph
$\Gamma_c$ on which the quantum graph is built.  We prepare with two simple lemmas.

\newpage
\begin{lemma}\label{Sandtaulemma}
Suppose that $\mathcal{S} x = x$.  Then the following equalities hold:
 \begin{itemize}
 \item
 $\mathcal{S}^* x = x$;
 \item
 $\mathcal{S} \tau x = \tau x$; and
 \item
 $\mathcal{S}^* \tau x = \tau x$.
 \end{itemize}
\end{lemma}

\begin{proof}
 If $\mathcal{S} x = x$, then
 $\mathcal{S}^* \mathcal{S} x = \mathcal{S}^* x$, and since $\mathcal{S}$ is unitary,  $x = \mathcal{S}^* x$.

Since  $\tau \mathcal{S} \tau =  \mathcal{S}^*$, $x = \tau \mathcal{S} \tau x$ so $\tau x = \mathcal{S} \tau x$.

The final claim follows by a similar argument.
\end{proof}

By replacing $x$ with $x \pm \tau x$ we can therefore assume that the set of fixed vectors of $\mathcal{S}$ has a basis such that $\tau x = \pm x $.

\begin{lemma}\label{LoBehold}
Consider any $k$-cycle $\mathcal{C} \in \Gamma_c$ and assign vertices $v_\ell$ and edge orientations so that $\{{\bf e_1}, {\bf e_2}, \dots, {\bf e_k}\} =  \{v_1 \to v_2,  v_2 \to v_3,  \dots v_k \to v_1\}$.  Then the vector  $x \in \mathbb{C}^{2m}$ such that $x({\bf e}_\ell) = 0$ if  ${\bf e}_\ell \notin \mathcal{C}$;
$x({\bf e}_\ell) = 1$ for $\ell = 1, \dots k$; and $x({\bf e}_\ell) = -1$ for ${\bf e}_\ell = \tau {\bf e}_j$, $j = 1, \dots k$, is a fixed point of
$\mathcal{S}$.
\end{lemma}

\begin{proof}
Direct verification.
\end{proof}

\begin{thm}\label{FixedSvecs}
We suppose that $\Gamma_c$ is finite and connected.
\begin{enumerate}[(A)]
\item
If $\mathcal{S}x = x$ and $\tau x = x$, then $x$ is a multiple of the constant vector ${\bf 1}  = ( 1,1,1, \dots)$.
\item
If $\mathcal{S}x = x$ and $x \perp {\bf 1}$, then $\tau x = - x$.
\item
The fixed points of $\mathcal{S}$ are harmonic $1$-forms in the sense that
\begin{equation}\label{harm1form}
(\extd_t -\extd_s) \mathcal{D}^{-1}(\extd^*_t - \extd^*_s ) x = 0.
\end{equation}
\item
There are $\beta_1 := m-n+1$
linearly independent harmonic $1$-forms such that $\tau x = - x$.
\item
Let $\{\mathcal{C}_j: 1 \le j \le \beta_1\}$ be a set of independent generators of the fundamental group of $\Gamma_c$.   (These are cycles, or ``holes''  in  $\Gamma_c$.)  Then there is a basis for the fixed vectors of $\mathcal{S}$ such that $\tau x = - x$ consisting of
vectors supported on each $\mathcal{C}_j$, which can be normalized to have the explicit forms
given in
Lemma \ref{LoBehold}.

\end{enumerate}
\end{thm}

\begin{rem}
The operator on the left side of \eqref{harm1form}
is the normalized Laplacian on discrete $1$-forms, also known as the ``{\em 1-down Laplacian},'' and thus
such $x$ are
discrete harmonic $1$-forms, or, in the terminology of \cite{BeLo}, ``{\em circulations}.''  The operator $\extd_s\mathcal{D}^{-1}\extd^*_t $
appearing in \eqref{Sdef}
is an asymmetric variant of the normalized Laplacian on forms as given on the left side of \eqref{harm1form}.  Parts of Theorem \ref{FixedSvecs} is
standard lore about harmonic $1$-forms, which we have adapted
 to account for edge orientations and to connect these objects to $\mathcal{S}$.  See  \cite{Sun08,Lim20} for more about
 discrete harmonic $1$-forms.
\end{rem}

\begin{proof}

\begin{enumerate}[(A)]
\item
Suppose that $\mathcal{S}x = x$ and $\tau x = x$.   It follows that $\extd_s \mathcal{D}^{-1}\extd^*_t x =  x$.
We calculate
\begin{align*}
\|x\|^2 &= \left\langle \extd^*_s x , \mathcal{D}^{-1}\extd^*_t x  \right\rangle \\
&= \sum_{{\bf e},{\bf e^\prime}}  \overline{x_{\bf e^\prime}} \delta_{s({\bf e^\prime}),t(({\bf e})} \frac{1}{d_{t(({\bf e})}} x_{\bf e} \\
&= \sum_{v \in \mathcal{V}} \frac{1}{d_v} \sum_{{\bf e} \to v} \sum_{{\bf e^\prime} \leftarrow v} \overline{x_{\bf e^\prime}} x_{\bf e}\\
&= \sum_{v \in \mathcal{V}} \frac{1}{d_v} \left|\sum_{{\bf e} \to v} x_{\bf e}\right|^2
\end{align*}
By Cauchy's inequality,
\[
\left|\sum_{{\bf e} \to v} x_{\bf e}\right|^2
 \le
\sum_{{\bf e} \to v} 1  \sum_{{\bf e} \to v} |x_{\bf e}|^2
= d_v \sum_{{\bf e} \to v} |x_{\bf e}|^2,
\]
with equality iff all $x_{\bf e}$ are equal whenever $t({\bf e}) = v$ and whenever $s({\bf e}) = v$.  Because the graph is connected, if equality holds for all vertices in the sum, it follows that $x$ is proportional to ${\bf 1}$.  We thus have

\begin{align*}
\|x\|^2 &\le \sum_{v \in \mathcal{V}} \sum_{{\bf e} \to v} |x_{\bf e}|^2 \\
 &= \|x\|^2,
\end{align*}
implying that $x$ is proportional to ${\bf 1}$.

\item
Claim B follows immediately from Claim A, since the states that are symmetric with respect to edge reversal are unique up to normalization.
\item
$x \perp {\bf 1}$, because
\begin{align*}
0 &= \left\langle {\bf 1}, \left(\mathcal{S}+\tau \right) x\right\rangle\\
   &=  \left\langle \left(\mathcal{S}^*+\tau \right) {\bf 1} , x\right\rangle\\
   &= \left\langle 2 {\bf 1} , x\right\rangle.
\end{align*}
Since
$\extd_s\mathcal{D}^{-1}\extd^*_t x = 0$, by taking the inner product with $\tau x$, we get
\begin{align*}
0 &= \left\langle \tau x , \extd_s \mathcal{D}^{-1}\extd^*_t x \right\rangle\\
   &= \left\langle  \extd_s^* \tau x , \mathcal{D}^{-1}\extd^*_t x \right\rangle\\
   &= \left\langle  \extd_t^* x , \mathcal{D}^{-1}\extd^*_t x \right\rangle\\
   &= \sum_v \frac{1}{d_v} \left|\sum_{t({\bf e}) = v} x_{\bf e}\right|^2 \ge 0.
\end{align*}
It follows that  $\left|\sum_{t({\bf e}) = v} x_{\bf e}\right|^2 = 0$, and a similar argument gives
$\left|\sum_{s({\bf e}) = v} x_{\bf e}\right|^2 = 0$, from which \eqref{harm1form} follows.

\item
This is clear if $x \propto {\bf 1}$.
Since the fixed vectors of $\mathcal{S}$ orthogonal to ${\bf 1}$ are all odd with respect to $\tau$, it follows that for them,

\begin{equation}\label{harmforms}
\extd_s\mathcal{D}^{-1}\extd^*_t x = 0,
\end{equation}
and as in the previous proof, on the subspace of $\mathbb{C}^{2m}$ for which $\tau x = - x$, $\extd_s\mathcal{D}^{-1}\extd^*_t x$
implies \eqref{harm1form}.

\item

According to (C), the dimension of the null space of $\mathcal{S} - \mathbb{1}$ when restricted to the subspace of $x: \tau x = - x$
equals that of
$(\extd_t -\extd_s)\mathcal{D}^{-1}(\extd^*_t -\extd^*_s)x$,
and because the rank of the Laplacian on discrete $1$-forms equals that of the scalar discrete Laplacian, which has nullity $1$ when the graph is connected, the
set of solutions of \eqref{harmforms} such that $\tau x = - x$
has dimension
$\beta_1$, c.f. for example \cite{FKW}, Theorem 20.
\item
Because of Lemma \ref{LoBehold}
we may choose the cycles on which the basis vectors are supported as any set of generators of the fundamental group of $\Gamma_c$.
\end{enumerate}

\end{proof}

\bigskip\bigskip\bigskip
\section{Topological quantum-graph states}\label{s: TQGS}

Due to its relative simplicity Kottos and Smilansky concentrated on a version of
the secular equation in the form
\eqref{SecEq1}
based on the bond-scattering matrix in \cite{KS99}, but they also derived an alternative ``vertex-scattering'' variant.
Their $n \times n$ vertex-scattering secular matrix is smaller than $2m \times 2m$, but depends on the eigenparameter $k$ in a more complicated way.  Moreover, rather than being an analytic function of $k$, it has singularities.
In Eq. (8) of \cite{KS99}, Kottos and Smilansky discarded these singularities with the consequence that potential quantum-graph states
may be overlooked
if their eigenvalues happen to coincide with the singularities.
 In this section we show that a factorization used by Ihara in \cite{Iha} provides a direct way to pass from one version of the secular matrix to the other, and
makes a deep connection with these possible eigenstates on the quantum graph, which are
linkd to its topological structure.
We turn our attention to the vertex-scattering secular equation below, cf. \eqref{VertSec}.

Because of Proposition \ref{Sdef}, in the usual ``bond-scattering'' form the secular determinant can be written
as
\begin{align}
\det (\mathbb{1} - \Phi \hat{\mathcal{S}}) &=
\det(\mathbb{1}_{4m} - \Phi (2 \extd_s \mathcal{D}^{-1}\extd_t^* - \tau ) \otimes \mathbb{1}_2)\nonumber\\
&= \det (\mathbb{1}_{4m}+ \Phi \hat{\tau} - 2 \Phi \extd_s \mathcal{D}^{-1}\extd_t^*\otimes \mathbb{1}_2 ).
\end{align}
(As in \eqref{SecEq}, $\hat{\mathcal{S}} := \mathcal{S} \otimes \mathbb{1}_{2}$.)
As we spell out below, an ``Ihara-style'' factorization reducing a
$4m \times 4m$ determinant to an $2n \times 2 n$ determinant
can be carried out if $\det(\mathbb{1} + \Phi \hat{\tau})$ can be treated as an invertible factor.  Since, however,
$(\mathbb{1} + \Phi \tau)$ is
not always invertible,
treating it as a factor to be inverted introduces singularities.
Generically this is not problematic, since for a generic quantum graph it is the zeros of the vertex-version of the secular determinant that yield eigenvalues of the quantum graph, and zeros cannot coincide with singularities. However, in non-generic cases,
the null space of
$\mathbb{1} + \Phi \tau$ and that of $\mathbb{1} - \Phi \hat{\mathcal{S}}$ may have a nontrivial intersection,
in which case the introduced singularities can coincide with eigenvalues of the quantum graph.
When this occurs, analytic aspects of the spectrum are entangled with the combinatorial topology of the graph, and hence the corresponding eigenvectors and eigenfunctions
can be viewed as {\em topological quantum-graph states}, cf. Def \ref{TQGSDeff}, below.

We capture some useful properties of
$\hat{\mathcal{S}}$ and $\Phi$ in the following lemma:

\begin{lemma}\label{Basicprops}
Let $x\in\mathbb{C}^{4m}$, thought of as a set of
2-vectors attached to oriented edges.  Then any two of the following imply the third:
\begin{enumerate}
\item $\hat{\mathcal{S}} x = x$,
\item $\Phi(k) x = x$,
\item $\Phi(k) \hat{\mathcal{S}} x = x$.
\end{enumerate}
\smallskip

\noindent
Moreover,
\begin{equation}\label{InverseCalc}
(\mathbb{1} + \Phi \tau)(\mathbb{1} - \Phi \tau) =
(\mathbb{1} - \Phi \tau)(\mathbb{1} + \Phi \tau) = \mathbb{1}  - \Phi \tau \Phi \tau,
\end{equation}
and hence any vector $x$ in the null space of $(\mathbb{1} \pm \Phi \tau)$
is a fixed vector of the block-diagonal matrix with blocks $\Phi_{{\bf e}}(k) \Phi_{-\bf e}(k)$, i.e.,
\begin{equation}\label{PerSols}
\Phi_{{\bf e}}(k) \Phi_{-\bf e}(k) x({\bf e}) = x({\bf e}).
\end{equation}
Eq.
\eqref{InverseCalc} provides an explicit formula for the inverse of $\mathbb{1} + \Phi \tau$, when
there are no non-trivial solutions to \eqref{PerSols}:
\begin{align}\label{e:inverse}
(\mathbb{1} + \Phi \tau)^{-1} &= \left(\mathbb{1}  - \Phi \tau\right) \left(\mathbb{1}  - \Phi \tau \Phi \tau\right)^{-1}\nonumber \\
&=  \left(\mathbb{1}  - \Phi \tau\right) \left(\mathbb{1}  -
\textrm{BMat}[\Phi_{\bf e}(k) \Phi_{-{\bf e}}(k)] \right)^{-1},
\end{align}

\end{lemma}

\begin{rem}
The expression \eqref{PerSols} shows that any vector in
the null space of $(\mathbb{1} \pm \Phi \tau)$ corresponds to a
periodic solution of the underlying ordinary differential equation on the doubled edge
obtained by following $-{\bf e}$ by ${\bf e}$.  We refer to these as {\em periodic
solutions on the doubled edge}, and caution that they may be nonzero on certain edges
while vanishing on others.  Similarly, the associated values of $k^2$ will be called
{\em doubled-edge periodic eigenvalues}.
In works where $V(x) = 0$ (e.g., \cite{BeKu})
these periodic solutions show up
as factors of the form $\exp(2 i k \ell_{\bf e})$.

\end{rem}

\begin{proof}
The first three statements and \eqref{InverseCalc} are immediate.
From \eqref{InverseCalc}
and the definition of $\tau$ the matrix $\Phi \tau \Phi \tau$ consists of blocks of the form Eq.
\eqref{PerSols}.
The calculation of the inverse is straightforward.
\end{proof}

We shall show in Corollary \ref{c:v-secular-m} that vectors in the null space of
$(\mathbb{1} + \Phi \hat{\tau})$ give rise to singularities in the vertex secular equation.  Since
eigenstates of the quantum graph are determined through their
counterparts $x$ on the associated combinatorial graph via
 \begin{equation*}
(\mathbb{1} - \Phi \hat{\mathcal{S}}) x = 0,
 \end{equation*}
 if they are also in this null space, i.e.,
 \begin{equation*}
 (\mathbb{1} + \Phi \hat{\tau})x = 0,
  \end{equation*}
 then the vertex secular equation fails to account for them.
Note that the difference of the these two conditions  yields
$ \Phi(k) \left(\hat{\mathcal{S}} + \hat{\tau}\right) x = 0$,
 which, since $\Phi(k)$ is invertible,
 is equivalent to the purely topological condition
\begin{equation}\label{topcon}
\left(\hat{\mathcal{S}} + \hat{\tau}\right) x = 0.
\end{equation}
This allows us to identify the special eigenstates that satisfy an additional purely topological condition:

%

\begin{defi}\label{TQGSDeff}
A quantum-graph eigenfunction $\psi$ such that its counterpart $x \in \mathbb{C}^{4m}$
simultaneously satisfies
 $\left(\hat{\mathcal{S}} + \hat{\tau}\right) x = 0$
and
 $(\mathbb{1} - \Phi \hat{\mathcal{S}}) x = 0$
will be called a {\em topological quantum-graph state}.
(Recall that by convention, on each directed edge ${\bf e}$,
$x_1({\bf e}) = \psi(s({\bf e}))$ and $x_2({\bf e}) = \psi^{\prime}(s({\bf e}))$.
\end{defi}

\begin{thm}\label{t:topoCharacterization}
Suppose that $\psi$ is a topological quantum-graph state and that
$x \in \mathbb{C}^{4m}$ is its counterpart edge-space vector.
Then
\begin{enumerate}[(A)]
\item\label{St1}

For $ \ell = 1,2$,
\begin{equation}\label{diverg}
\extd_t^* x_\ell = 0.
\end{equation}
\item\label{St2}
On every edge, $x_1({\bf e}) = 0$.
Consequently,
$\psi$ satisfies Dirichlet conditions at every vertex.  (This statement
excludes ``fake vertices'' of degree 2.)
In particular, topological quantum-graph states vanish on {\em leaves}, which are by definition edges having a terminal degree-one vertex.
\item\label{DBConCycles}
There are at most $\beta_1 = m-n+1$
topological quantum-graph states for a given eigenvalue, and they have an eigenbasis consisting of
vectors supported on generating cycles $\mathcal{C}_j$ of the fundamental group.
For any topological quantum-graph state in such a basis, the Dirichlet problems for the edges ${\bf e} \subset \mathcal{C}_k$
on which it is supported
must share the eigenvalue $k^2$.
\end{enumerate}
\end{thm}

\begin{rems}
\,
\begin{itemize}
\item
The expression on the left side of \eqref{diverg} can be thought of as a discrete divergence.
\item
The final statement of
\ref{DBConCycles} means that even if the quantum graph has a nontrivial interaction, the
edges comprised in $ \mathcal{C}_k$ are
``spectrally equilateral.''
\item
In the case of equilateral quantum graphs the dimension of the eigenbasis is equal to $\beta_1$.
\end{itemize}
\end{rems}

\bigskip

 \begin{proof}\,

 \par

 \ref{St1}:
It  follows from Definition \ref{TQGSDeff} together with Eq. \eqref{Sdef} that
\begin{equation}\label{twisted}
\extd_s \mathcal{D}^{-1} \extd_t^* x_\ell = 0.
\end{equation}
Taking the inner product of $\extd_s \mathcal{D}^{-1}\extd_t^* x_\ell$ with $\tau x_\ell$ yields
\begin{align}
0 &= \left\langle{\tau x_\ell, \extd_s \mathcal{D}^{-1}\extd_t^* x_\ell}\right\rangle\nonumber\\
&= \left\langle{\extd_t^* x_\ell,  \mathcal{D}^{-1}\extd_t^* x_\ell}\right\rangle,
\end{align}
which implies that $\extd_t^* x_\ell = 0$, because $\mathcal{D}^{-1} > 0$.  (Note:  For $\ell=2$ \eqref{twisted} is also just
the Kirchhoff condition on the derivatives at the vertices.)

\ref{St2}:
Since
\begin{align*}
[\extd_t^* x_1](v) &=  \sum_{{\bf e}: v = t({\bf e})} {x_1({\bf e})} \\
&= d_v \psi(v),
\end{align*}
It follows from Statement \ref{St1} that $\psi(v) = 0$ at every vertex.
The final statement then follows since at a terminal vertex, $\psi(v) = \psi^\prime(v) = 0$ and hence $x = 0$ on the associated edge.

\ref{DBConCycles}:
It follows from Statement \ref{St1} that
\begin{equation}\label{harm1formhat}
(\extd_t -\extd_s) \mathcal{D}^{-1}(\extd^*_t - \extd^*_s ) x_2 = 0,
\end{equation}
so $x_2$ is a harmonic $1$-form in the same sense as
in  \eqref{harm1form}.  Since the Laplacian as written in self-adjoint form in \eqref{harm1formhat}, its null space has an orthonormal basis.  The constant vector
in $\mathbb{C}^{2m}$ is in the null space but cannot be a topological quantum graph space.  Hence the dimension of the null space is determined as in the proof of
Theorem \ref{FixedSvecs}.  The statement about basis vectors supported on the generators of the fundamental group
follows also as for Theorem \ref{FixedSvecs} because any topological quantum-graph state must be nontrivial on some cycle, and linear combinations can always be found
with support on simple cycles.
 \end{proof}

\begin{rem}
If the edge interactions are symmetric under the action of $\tau$, the counterparts $x$ of the basis vectors on the cycles
are of the form given in Lemma \ref{LoBehold}.
Their existence of course requires the existence of cycles consisting of sequences of edges with a common Dirichlet eigenvalues.
Expressed in the more general schema of the operators $\phi_{\bf e}$, topological quantum-graph states require that each $\phi_{\bf e}$
for ${\bf e}$ in the cycle share the
eigenvector
$
\begin{bmatrix}
    0 \\
   1
\end{bmatrix}
$,
and if the supporting cycle consists of the consecutive oriented edges ${\bf e}_1 \dots {\bf e}_1$, then
\[
\phi_{{\bf e}_n} \phi_{{\bf e}_{n-1}}  \dots \phi_{{\bf e}_1}
\begin{bmatrix}
    0 \\
   1
\end{bmatrix}
=
\begin{bmatrix}
    0 \\
   1
\end{bmatrix}.
\]
\end{rem}

\begin{thm} \label{t:classification}
There are three types of eigenvalues
of a quantum graph:
\begin{itemize}
\item  The eigenvalues of topological eigenstates are periodic eigenvalues associated with doubled edges
in accordance with \eqref{PerSols}.  However, not all  periodic eigenvalues associated with doubled edges
give rise to topological eigenstates, only those that are associated with spectrally equilateral cycles as in
Theorem \ref{t:topoCharacterization}, part \ref{DBConCycles}.
\item
Continuing to use the notation $\hat{M} := M \otimes \mathbb{1}_2\,$  for any matrix $M$,
the values $k$ for which the
{\em vertex secular matrix}
 \begin{equation}\label{e:v-scattering-secular}
\widehat{\extd_t^*} \left(  \textrm{BMat}(\Phi_{\bf e}(k) + \Phi_{-{\bf e}}(k)^{-1}) - 2  \hat{\tau}\right)
\textrm{BMat}(\Phi_{\bf e}(k) - \Phi_{-{\bf e}}(k)^{-1})
\widehat{\extd_t}
\end{equation}
has zero determinant
correspond to {\em non-topological eigenstates} with eigenvalue $k^2$.
\item
If $k^2$ is a periodic eigenvalue of one or more doubled edges,
it is possible for $k^2$ to simultaneously be in the topological spectrum and in the nontopological
spectrum, cf. Case Study \ref{AccDeg}, below.
\end{itemize}
\end{thm}

\begin{rem}
Observe that the components $\secmat_{vw}$
of the
{\em vertex secular matrix} $\secmat:= \widehat{\extd_t^*} \left( 2  (\mathbb{1} + \Phi \hat{\tau})^{-1}- \mathbb{1} \right)\widehat{\extd_t}
$ vanish unless vertices $v$ and $w$ are equal or connected by an edge.  Therefore
$\secmat$ is like a $k$-dependent discrete Laplacian matrix with weights that are not necessarily positive,
and which become singular for special values of $k$.
\end{rem}

\begin{proof} The characterization of the the topological eigenvalues and  eigenstates recapitulates part of Theorem \ref{t:topoCharacterization}.

To show the second part of the theorem,
which extends results  of \cite{KS99} for metric graphs where $V=0$,
suppose that $\det(\mathbb{1} + \Phi(k) \tau) \neq 0$.
Following the proof of Ihara's Theorem as presented in \cite{kempton2016non} we begin writing
\begin{align}\label{1stStep}
\det (\mathbb{1} - \Phi \hat{\mathcal{S}} ) &= \det(\mathbb{1} - \Phi (2 \extd_s \mathcal{D}^{-1}\extd_t^* - \tau )\otimes \mathbb{1}_2)\nonumber\\
&= \det (\mathbb{1} + \Phi \hat{\tau} - 2 \Phi (\extd_s  \otimes \mathbb{1}_2) (\mathcal{D}^{-1}\otimes \mathbb{1}_2) (\extd_t^* \otimes \mathbb{1}_2 ))
\nonumber\\
&= \det (\mathbb{1} + \Phi \hat{\tau} - 2 \Phi \widehat{\extd_s} \widehat{\mathcal{D}^{-1}} \widehat{\extd_t^*}).
\end{align}
Next we use the assumption that $\mathbb{1} + \Phi(k) \tau$ is invertible and the matrix identity
\begin{equation}\label{MatID}
\det(Q + MN) = \det(Q)\det(\mathbb{1} + NQ^{-1}M),
\end{equation}
the validity of which requires
only that the dimensions of $M, N, Q$, and $\mathbb{1}$ allow them to be added and multiplied as shown, plus
the assumption that $Q$ is square and invertible.
If $k^2$ is not
a doubled-edge periodic eigenvalue,
this enables us to rewrite \eqref{1stStep} as
\begin{align}\label{e:factorization}
\det (\mathbb{1} - \Phi \hat{\mathcal{S}} ) &= \det(\mathbb{1} + \Phi \hat{\tau}) \det(\mathbb{1} - 2 \widehat{\mathcal{D}^{-1}} \widehat{\extd_t^*} (\mathbb{1} + \Phi \hat{\tau})^{-1}  \Phi  \widehat{\extd_s})\nonumber\\
&= \det(\mathbb{1} + \Phi \hat{\tau}) \det(\mathbb{1} - 2 \widehat{\mathcal{D}^{-1}} \widehat{\extd_t^*} (\mathbb{1} + \Phi \hat{\tau})^{-1}  \Phi \hat{\tau} \widehat{\extd_t})\nonumber\\
&= \det(\mathbb{1} + \Phi \hat{\tau}) \det\left(\widehat{\mathcal{D}^{-1}}\widehat{\extd_t^*} \left(\mathbb{1} - 2  (\mathbb{1} + \Phi \hat{\tau})^{-1}  \Phi  \hat{\tau} \right)\widehat{\extd_t}\right)\nonumber\\
&= \det(\mathbb{1} + \Phi \hat{\tau}) \left(\det \mathcal{D}\right)^{-2} \det \left(\widehat{\extd_t^*} \left(\mathbb{1} - 2  (\mathbb{1} + \Phi \hat{\tau})^{-1}  \Phi  \hat{\tau}\right)\widehat{\extd_t}\right)\nonumber\\
&= \det(\mathbb{1} + \Phi \hat{\tau}) \left(\det \mathcal{D}\right)^{-2} \det \left(\widehat{\extd_t^*} \left( 2  (\mathbb{1} + \Phi \hat{\tau})^{-1}- \mathbb{1} \right)\widehat{\extd_t}\right).
\end{align}
(The first equality uses \eqref{MatID}, the next two lines follow from identities in Theorem \ref{Facto}, and the final lines
are elementary.)
Thus the solutions in $k$ of
 \begin{equation}
  \det \left(\widehat{\extd_t^*} \left( 2  (\mathbb{1} + \Phi \hat{\tau})^{-1}- \mathbb{1} \right)\widehat{\extd_t}\right) = 0
\end{equation}
define non-topological eigenstates with eigenvalue $k^2$.

If $\det(\mathbb{1} + \Phi(k) \tau) = 0$, then either the
vertex secular equation is singular and cannot simultaneously be zero, or else it must contain at least one root associated with an eigenvalue,
in order to cancel the singularity.

The third type of eigenvalue arises due to a similar observation.  Since the edge secular equation is entire and the vertex secular equation is equivalent to dividing it by $\det (\mathbb{1} + \Phi(k) \tau)$, the zeroes of which correspond to periodic eigenvalues on the doubled edges, the degree
of a pole of the vertex secular equation is the difference of the degrees of the zeroes of $\det (\mathbb{1} + \Phi(k) \tau)$ and of the
bond secular equation.

\end{proof}

\begin{cor}\label{c:v-secular-m}
The matrix
\begin{equation}\label{e:v-secular-m}
\mathfrak{A}: =\widehat{\extd_t^*} \left( 2  (\mathbb{1} + \Phi \hat{\tau})^{-1}- \mathbb{1} \right)\widehat{\extd_t}
\end{equation}
 has the form
\begin{equation}
\mathfrak{A}_{vw} = \begin{cases}
2 \sum_{\{\textbf{e}: v \rightarrow w\}}\, (\Phi_{-\mathbf{e}}^{-1}- \Phi_{\mathbf{e}})^{-1} \text{, if } v \neq w \\
\sum_{\{\textbf{e}: \textbf{e} \rightarrow v\}}(\Phi_{-\mathbf{e}}^{-1}+ \Phi_{\mathbf{e}})(\Phi_{-\mathbf{e}}^{-1}- \Phi_{\mathbf{e}}) \text{, if } v= w.
\end{cases}
\end{equation}
\end{cor}

\begin{proof}
We recall the formula \eqref{e:inverse} and we take a look at the expression between the $\extd_t^*$ and $\extd_t$:
\begin{multline}
 2  (\mathbb{1} + \Phi \hat{\tau})^{-1}- \mathbb{1}
 = (2(\mathbb{1}- \Phi \hat{\tau}) - \text{BMat}[\mathbb{1}- \Phi_{\mathbf{e}}\Phi_{-\mathbf{e}}]
 )(\text{BMat}[\mathbb{1}- \Phi_{\mathbf{e}}\Phi_{-\mathbf{e}}])^{-1}\\
 =
 (\mathbb{1} + \text{BMat}[\Phi_{\mathbf{e}}\Phi_{-\mathbf{e}}] - 2\Phi \hat \tau)(\text{BMat}[\mathbb{1}- \Phi_{\mathbf{e}}\Phi_{-\mathbf{e}}])^{-1}.
\end{multline}
Now \begin{equation}
(\mathbb{1}+ \Phi_{\mathbf{e}}\Phi_{-\mathbf{e}})(\mathbb{1}- \Phi_{\mathbf{e}}\Phi_{-\mathbf{e}})^{-1} = (\Phi_{-\mathbf{e}}^{-1}+ \Phi_{\mathbf{e}})(\Phi_{-\mathbf{e}}^{-1}- \Phi_{\mathbf{e}})\\
\end{equation}
and
\[
\Phi \hat \tau = \hat \tau(\hat \tau \Phi \hat \tau) = \hat \tau \text{BMat}(\Phi_{-\mathbf{e}}).
\]
Thus

\begin{multline}
\Phi \hat \tau (\text{BMat}[\mathbb{1}- \Phi_{\mathbf{e}}\Phi_{-\mathbf{e}}])^{-1} = \hat \tau \text{BMat}(\Phi_{-\mathbf{e}})(\text{BMat}[\mathbb{1}- \Phi_{\mathbf{e}}\Phi_{-\mathbf{e}}])^{-1}
\\
= \hat \tau (\text{BMat}[\Phi_{-\mathbf{e}}^{-1}- \Phi_{\mathbf{e}}])^{-1}.\\
\end{multline}
This yields that
\begin{equation}\label{e:1}
2  (\mathbb{1} + \Phi \hat{\tau})^{-1}- \mathbb{1} = \text{BMat}[(\Phi_{-\mathbf{e}}^{-1}+ \Phi_{\mathbf{e}})(\Phi_{-\mathbf{e}}^{-1}- \Phi_{\mathbf{e}})] - 2\hat \tau (\text{BMat}[\Phi_{-\mathbf{e}}^{-1}- \Phi_{\mathbf{e}}])^{-1}.
\end{equation}
Since an edge pointing at vertex $v$ cannot also point to vertex $w\neq v$, the term in \eqref{e:v-secular-m} resulting from the first term above is again a block diagonal matrix:
\begin{equation}
(\widehat{\extd_t^*} \text{BMat}[(\Phi_{-\mathbf{e}}^{-1}+ \Phi_{\mathbf{e}})(\Phi_{-\mathbf{e}}^{-1}- \Phi_{\mathbf{e}})]\widehat{\extd_t})_{vw}= \begin{cases}
\sum_{\{\textbf{e}: \textbf{e} \rightarrow v\}}(\Phi_{-\mathbf{e}}^{-1}+ \Phi_{\mathbf{e}})(\Phi_{-\mathbf{e}}^{-1}- \Phi_{\mathbf{e}}) \text{, if } v= w \\
0 \text{ otherwise}.
\end{cases}
\end{equation}
The term resulting from the second term in \eqref{e:1} is not a diagonal block matrix, since the $\hat \tau$ reverses the edge orientations. Thus this term will be as follows:
\begin{equation}
2(\widehat{\extd_t^*} \hat \tau (\text{BMat}[\Phi_{-\mathbf{e}}^{-1}- \Phi_{\mathbf{e}}])^{-1}\widehat{\extd_t})_{vw}= \begin{cases}
2 \sum_{\{\textbf{e}: v \rightarrow w\}}\, (\Phi_{-\mathbf{e}}^{-1}- \Phi_{\mathbf{e}})^{-1} \text{, if } v \neq w  \\
0 \text{ otherwise}.
\end{cases}
\end{equation}

\end{proof}

\begin{cor}
Suppose $V=0$. Then the non-topological eigenvalues of the quantum graph are $k \in \mathbb{R}$ such that  $\det \mathfrak{A}(k)=0$ with
\begin{equation}\label{VertSec}
\mathfrak{A}(k)_{v, w}
= i\begin{cases} 0 \text{ if } v \text{ is not connected to } w
\\
- \frac{1}{\sqrt{\deg v \deg w} \,\, \sin(kl_{v, w})} \text{ if } v \text{ is connected to } w \text{ by an edge of length } l_{v, w}
\\
 \sum_{y \sim v}\frac{\cot (kl_{v, y})}{\deg (v)} \text{ if } v = w.
\end{cases}
\end{equation}

\end{cor}

\begin{rem}
In this case, the matrix $\mathfrak{A}$ is the vertex-scattering matrix from Kottos-Smilansky.
\end{rem}

\begin{proof}

This follows directly from Corollary \ref{c:v-secular-m} and identification \eqref{e:i}.
\end{proof}

\begin{cor}
A quantum graph with $V=0$ and edge lengths which are not rationally related has only non-topological eigenvalues.
For any graph that contains a cycle with rationally related lengths, a topological eigenvalue and eigenvector can be constructed.
\end{cor}

\section{Case studies}

In this section we collect some illustrative examples, many of which are based on tetrahedra, which are accessible without being trivial.

\begin{enumerate}
\item\label{Tet}
The equilateral tetrahedron with edge length $\pi$ and no potential.  Every quadrilateral cycle supports topological quantum-graph states as in
Lemma \ref{LoBehold}, with eigenvalues $\ell^2$, $\ell$ any positive integer.  Every triangular face supports topological quantum-graph states
with eigenvalues $\ell^2$, $\ell >0$ even.  These are not all linearly independent, as the dimension of the eigenspaces is only
$\beta_1 =3$, as can either be
shown directly or deduced as a consequence of Theorems \ref{FixedSvecs} and \ref{t:topoCharacterization}.  The nontopological eigenfunctions correspond to
eigenfunctions that
are even with respect to permutation of the edges incident to any given vertex, and thus to functions on edges proportional to $\cos \ell x$, with
$\ell \ge 0$ even.

\begin{figure}[h]\label{f:tetrahedra}
\centering{\includegraphics[scale=0.25]{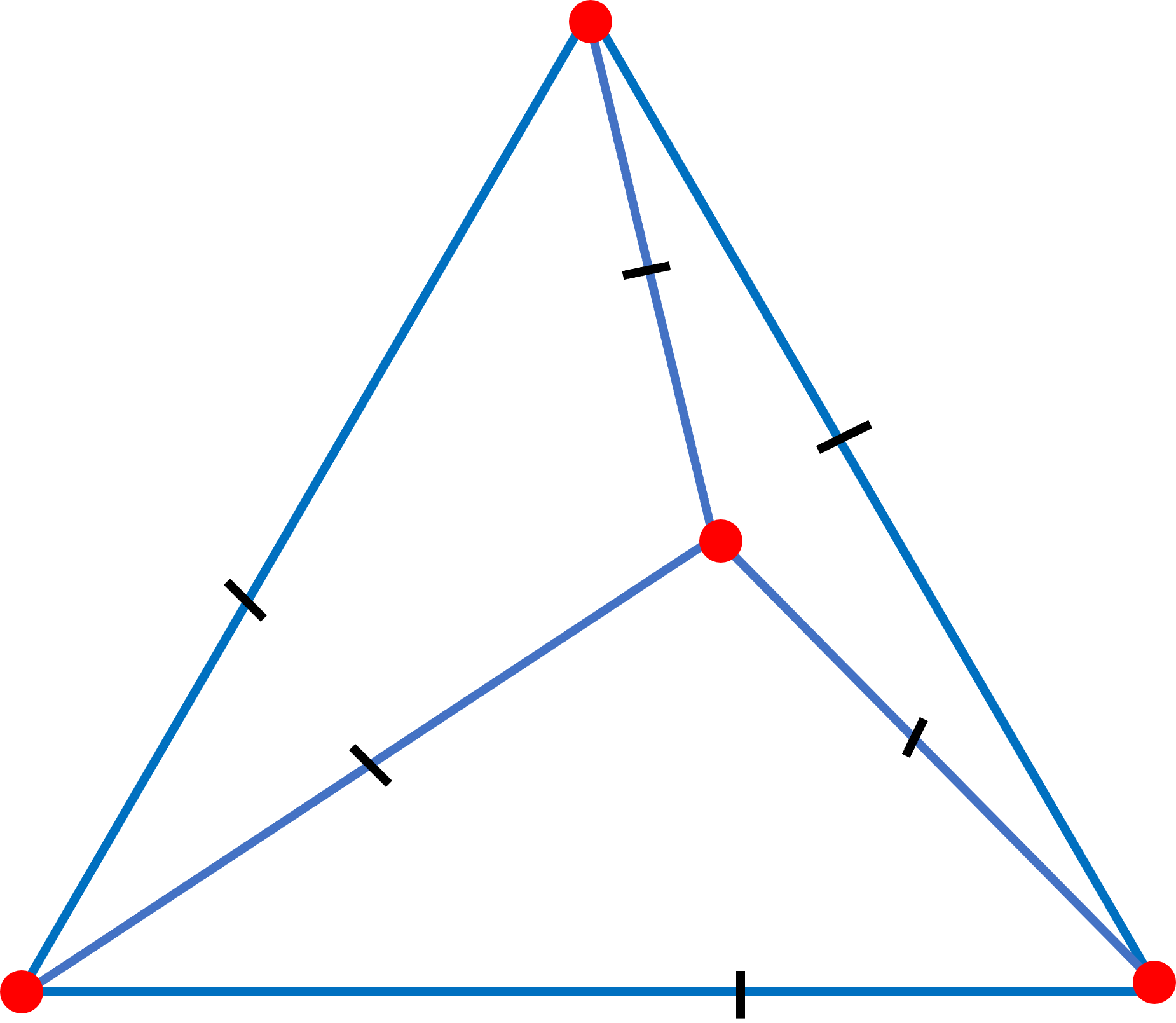}\includegraphics[scale=0.4]{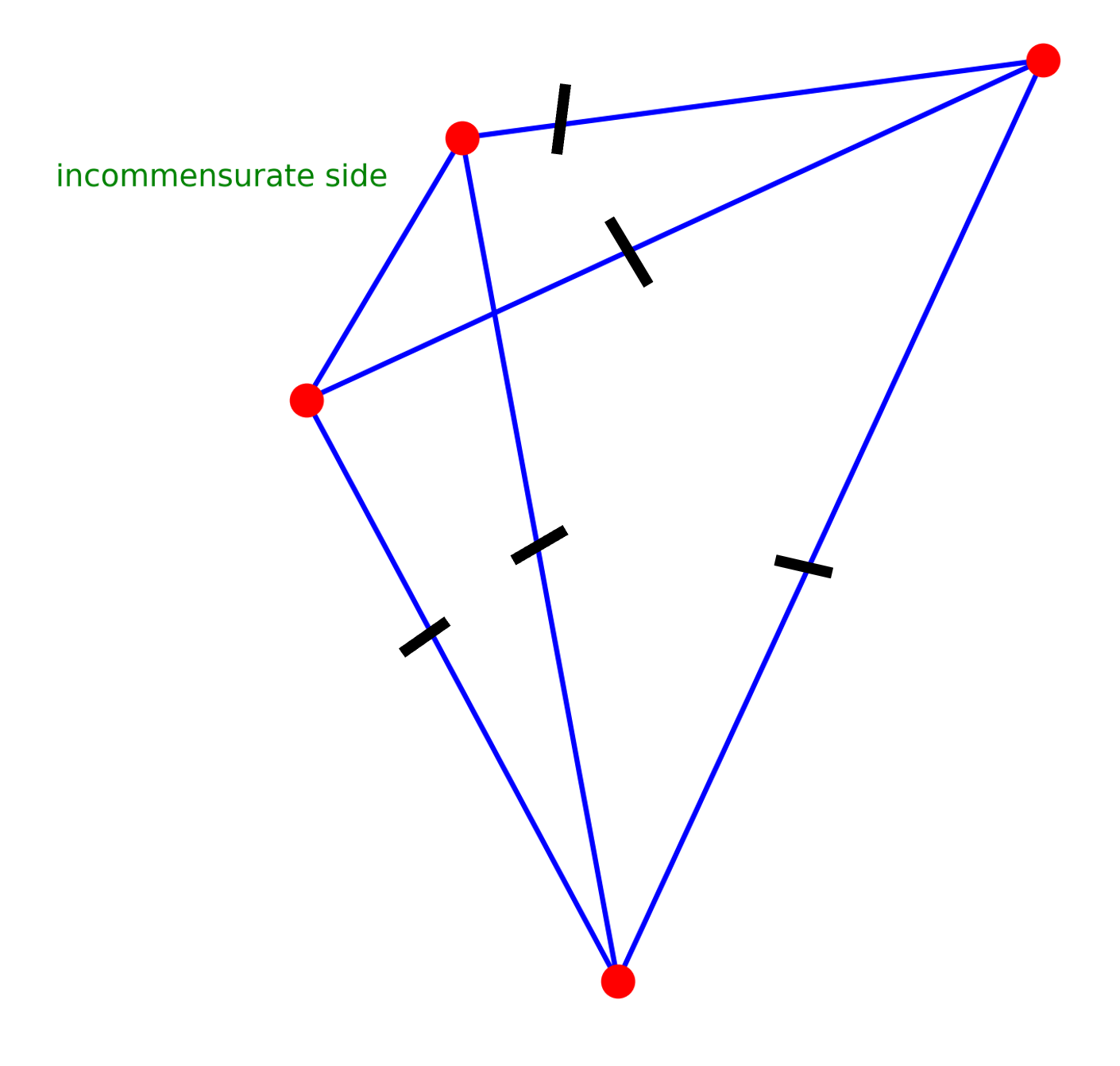}\includegraphics[scale=0.25]{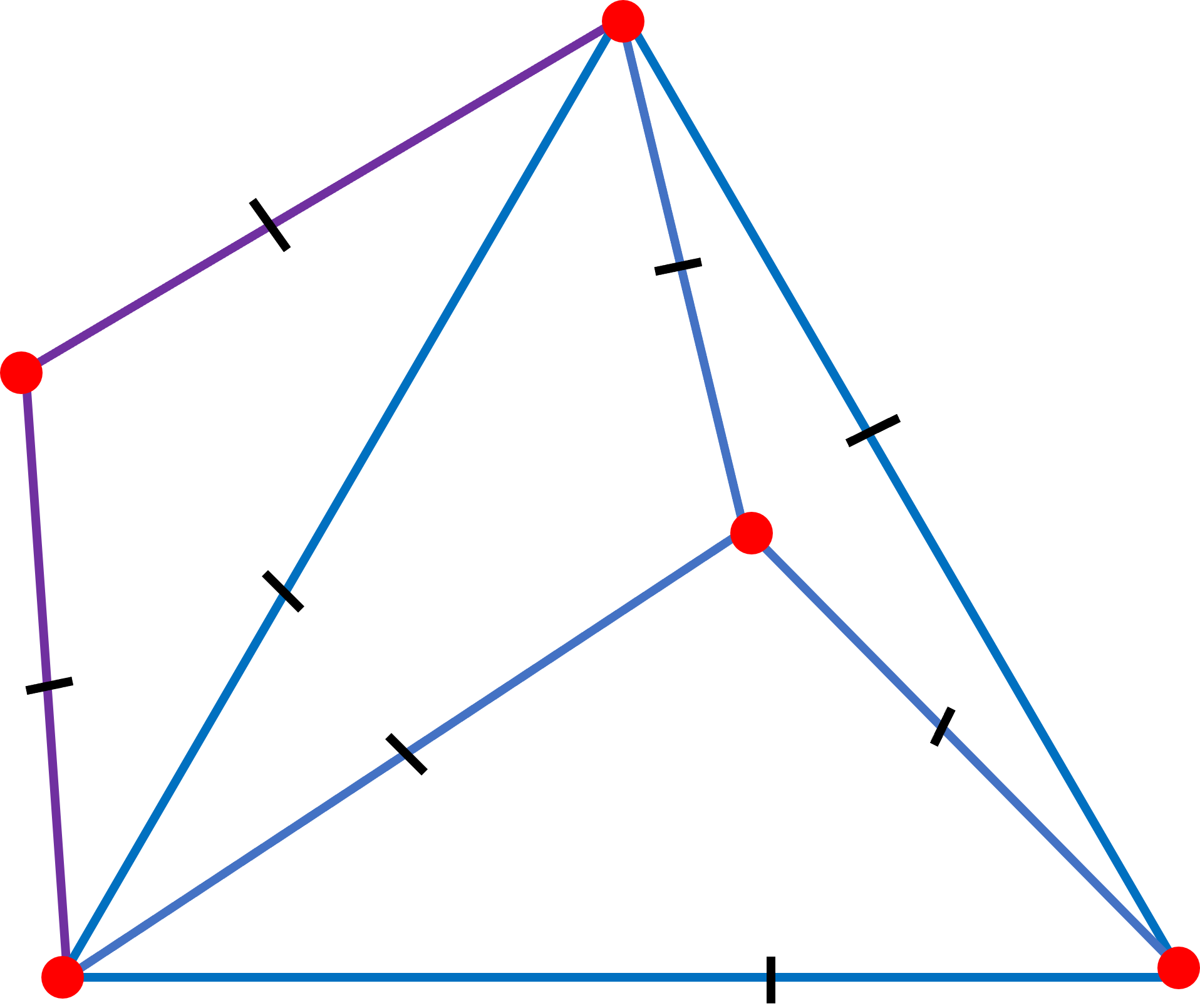}
\caption{Equilateral tetrahedron Case Study 1 (Left), tetrahedron with one incommensurate side from Case Study 2 (Center), and tetrahedron with two additional edges Case Study 3 (Right)}
}
\end{figure}

\item
A tetrahedral metric graph with edge lengths $\pi$ except for edge \#1, which is incommensurate, e.g., with length 1.  The topological states of the equilateral tetrahedron that are zero on edge \#1 persist, and they have multiplicity 2 rather than 3.
\item
If two additional edges of length $\pi$
are attached to vertices of the tetrahedron of Case Study 1, and their other ends are joined
at a new vertex of degree 2, we get a model that can be  thought of as equilateral, or it can be thought of as containing an additional edge of length $2 \pi$ joining the vertices.  All of the original topological states persist, and additional topological states appear that are
supported in part on the new edges.   The multiplicity is increased to the new $\beta_1 = 4$.
(The vertices to which the new edge is attached could be the same, in effect attaching a loop of lenght $2 \pi.$)  This confirms that the tetrahedral symmetry of Case Study 1 is not the source of the topological states.  Rather, they arise from equilaterality, as suitably interpreted.
\item
An equilateral length-$\pi$ tetrahedron, with $V(x) = 144$ on edge \#1, 0 on all others.  There are a full set of TQGSs supported on the cycles that do not include edge \#1, while on the cycles including edge \#1 there are only TQGSs corresponding to Pythagorean triples with the integer 12, i.e., (5,12,13), (9, 12, 15), (12, 16, 20).  This confirms that spectral equilaterality matters, not strict metric equilaterality.
\item
Although topological quantum-graph states can be regarded as probability-zero events in many randomly generated families of quantum graphs, cf. \cite{KS99},  there are other random models in which many of them exist with high probability. For instance, consider the model of
a large complete metric graph, for which the edge lengths are assigned one of two values according to some Bernoulli process. With a high probability the result will contain equilateral cycles as subgraphs, and these always support topological states according to Lemma \ref{LoBehold}. Theorem \ref{FixedSvecs} (D) implies that there will be a basis of topological bound states of dimension $\beta_1$. As the two-length model is a union of two complementary Erd\H{o}s-R\'enyi graphs, well-established work on estimates on Betti numbers gives the total number of edges in the graph as an asymptotic estimate on dimension of the eigenspace corresponding to the topological eigenvalue \cite{Kah}. Furthermore, by Theorem \ref{t:topoCharacterization}
a similar interacting model can be
made by assigning one of a pair of possible transfer matrices to the edges.
\begin{figure}[h]\label{f:two-lengths}
\centering{\includegraphics[scale=0.4]{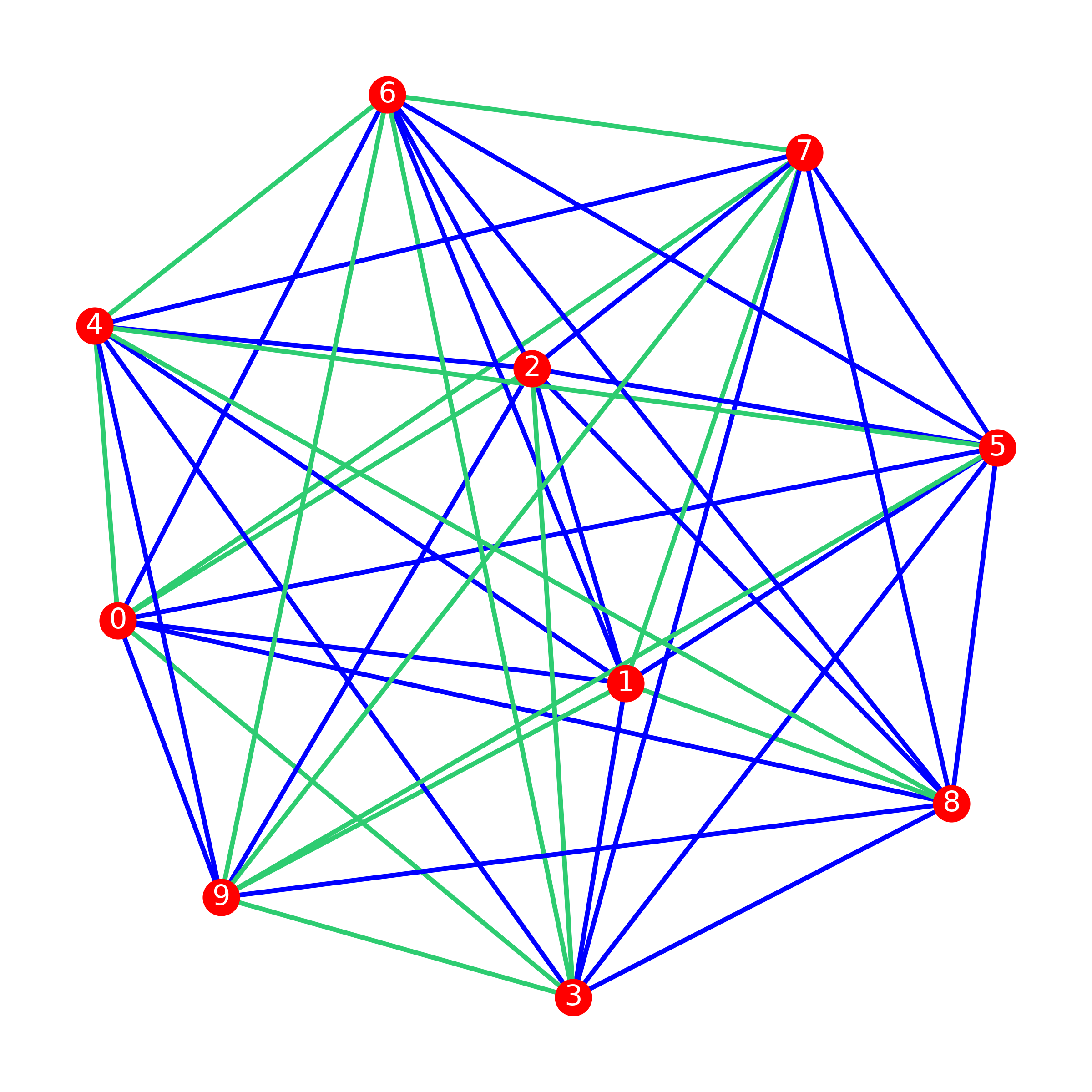}
\caption{A manifestation of the two-lengths model on 10 vertices where one length is marked in blue and the other in green, the green one chosen with probability 1/3}
}
\end{figure}

\item\label{AccDeg}
As stated in the third case of Theorem \ref{e:v-scattering-secular} a given eigenvalue may easily have both topological and ordinary eigenvectors.  Take, for instance, any quantum graph $\Gamma_1$ with non-topological eigenfunction $\psi$ that is not the ground state, with eigenvalue $k^2$.
Separately, consider an equilateral complete graph, such as the triangle or tetrahedron, with no potential energy, and adjust the scale so that one of its topological quantum-graph eigenvalues equals $k^2$.  To be specific, if the length of the edges is $L$, then $\frac{\pi^2}{L^2}$ is the eigenvalue of the topological states with eigenfunctions of the form $\sin(\pi \frac x L )$ on the edges of a cycle, so choose $L = \frac k \pi$.  Call this metric graph $T$.  Since $\psi$ is not the ground state of $\Gamma_1$, it vanishes somewhere, which we can designate as a degree-two vertex $v$.  Now ``surgically'' attach $v$ to one of the vertices of $T$ to form a new quantum graph $\Gamma_2$ for which both $\Gamma_1$ and $T$ can be regarded as subgraphs
(see \cite{BKKM} for a survey and analysis of surgery on quantum graphs.)  The original $\psi$ supported on $\Gamma_1$
continues to satisfy all of the conditions to be an eigenfunction on $\Gamma_2$ with eigenvalue $k^2$, and the same is true for the topological state supported on a cycle of $T$.
\end{enumerate}

\bigskip\bigskip

\begin{ack}
The authors thank Ginestra Bianconi for an extremely useful discussion.
In addition we are grateful to the referee for insights into topological resonances and related references.
\end{ack}

\begin{funding}
This work was partially supported by
the Royal Society University Research Fellowship  (grant
numbers RF$\backslash$ERE$\backslash$21005).
\end{funding}

\bibliographystyle{amsplain}

\end{document}